\documentclass[12pt]{amsproc}
\usepackage{amssymb}
\begin{document}
\theoremstyle{plain}
\newtheorem{thm}{Theorem}[section]
\newtheorem{theorem}[thm]{Theorem}
\newtheorem{lemma}[thm]{Lemma}
\newtheorem{corollary}[thm]{Corollary}
\newtheorem{proposition}[thm]{Proposition}
\newtheorem{addendum}[thm]{Addendum}
\newtheorem{variant}[thm]{Variant}
\theoremstyle{definition}
\newtheorem{notations}[thm]{Notations}
\newtheorem{question}[thm]{Question}
\newtheorem{problem}[thm]{Problem}
\newtheorem{remark}[thm]{Remark}
\newtheorem{remarks}[thm]{Remarks}
\newtheorem{definition}[thm]{Definition}
\newtheorem{claim}[thm]{Claim}
\newtheorem{assumption}[thm]{Assumption}
\newtheorem{assumptions}[thm]{Assumptions}
\newtheorem{properties}[thm]{Properties}
\newtheorem{example}[thm]{Example}
\numberwithin{equation}{thm}
\catcode`\@=11
\def\opn#1#2{\def#1{\mathop{\kern0pt\fam0#2}\nolimits}}
\def\bold#1{{\bf #1}}%
\def\underrightarrow{\mathpalette\underrightarrow@}
\def\underrightarrow@#1#2{\vtop{\ialign{$##$\cr
 \hfil#1#2\hfil\cr\noalign{\nointerlineskip}%
 #1{-}\mkern-6mu\cleaders\hbox{$#1\mkern-2mu{-}\mkern-2mu$}\hfill
 \mkern-6mu{\to}\cr}}}
\let\underarrow\underrightarrow
\def\underleftarrow{\mathpalette\underleftarrow@}
\def\underleftarrow@#1#2{\vtop{\ialign{$##$\cr
 \hfil#1#2\hfil\cr\noalign{\nointerlineskip}#1{\leftarrow}\mkern-6mu
 \cleaders\hbox{$#1\mkern-2mu{-}\mkern-2mu$}\hfill
 \mkern-6mu{-}\cr}}}
\let\amp@rs@nd@\relax
\newdimen\ex@
\ex@.2326ex
\newdimen\bigaw@
\newdimen\minaw@
\minaw@16.08739\ex@
\newdimen\minCDaw@
\minCDaw@2.5pc
\newif\ifCD@
\def\minCDarrowwidth#1{\minCDaw@#1}
\newenvironment{CD}{\@CD}{\@endCD}
\def\@CD{\def\A##1A##2A{\llap{$\vcenter{\hbox
 {$\scriptstyle##1$}}$}\Big\uparrow\rlap{$\vcenter{\hbox{%
$\scriptstyle##2$}}$}&&}%
\def\V##1V##2V{\llap{$\vcenter{\hbox
 {$\scriptstyle##1$}}$}\Big\downarrow\rlap{$\vcenter{\hbox{%
$\scriptstyle##2$}}$}&&}%
\def\={&\hskip.5em\mathrel
 {\vbox{\hrule width\minCDaw@\vskip3\ex@\hrule width
 \minCDaw@}}\hskip.5em&}%
\def\verteq{\Big\Vert&&}%
\def\noarr{&&}%
\def\vspace##1{\noalign{\vskip##1\relax}}\relax\let\amp@rs@nd@&\iffalse}\fi
 \CD@true\vcenter\bgroup\relax\let\\=\cr\iffalse}\fi\tabskip\z@skip\baselineskip20\ex@
 \lineskip3\ex@\lineskiplimit3\ex@\halign\bgroup
 &\hfill$\m@th##$\hfill\cr}
\def\@endCD{\cr\egroup\egroup}
\def\>#1>#2>{\amp@rs@nd@\setbox\z@\hbox{$\scriptstyle
 \;{#1}\;\;$}\setbox\@ne\hbox{$\scriptstyle\;{#2}\;\;$}\setbox\tw@
 \hbox{$#2$}\ifCD@
 \global\bigaw@\minCDaw@\else\global\bigaw@\minaw@\fi
 \ifdim\wd\z@>\bigaw@\global\bigaw@\wd\z@\fi
 \ifdim\wd\@ne>\bigaw@\global\bigaw@\wd\@ne\fi
 \ifCD@\hskip.5em\fi
 \ifdim\wd\tw@>\z@
 \mathrel{\mathop{\hbox to\bigaw@{\rightarrowfill}}\limits^{#1}_{#2}}\else
 \mathrel{\mathop{\hbox to\bigaw@{\rightarrowfill}}\limits^{#1}}\fi
 \ifCD@\hskip.5em\fi\amp@rs@nd@}
\def\<#1<#2<{\amp@rs@nd@\setbox\z@\hbox{$\scriptstyle
 \;\;{#1}\;$}\setbox\@ne\hbox{$\scriptstyle\;\;{#2}\;$}\setbox\tw@
 \hbox{$#2$}\ifCD@
 \global\bigaw@\minCDaw@\else\global\bigaw@\minaw@\fi
 \ifdim\wd\z@>\bigaw@\global\bigaw@\wd\z@\fi
 \ifdim\wd\@ne>\bigaw@\global\bigaw@\wd\@ne\fi
 \ifCD@\hskip.5em\fi
 \ifdim\wd\tw@>\z@
 \mathrel{\mathop{\hbox to\bigaw@{\leftarrowfill}}\limits^{#1}_{#2}}\else
 \mathrel{\mathop{\hbox to\bigaw@{\leftarrowfill}}\limits^{#1}}\fi
 \ifCD@\hskip.5em\fi\amp@rs@nd@}
\newenvironment{CDS}{\@CDS}{\@endCDS}
\def\@CDS{\def\A##1A##2A{\llap{$\vcenter{\hbox
 {$\scriptstyle##1$}}$}\Big\uparrow\rlap{$\vcenter{\hbox{%
$\scriptstyle##2$}}$}&}%
\def\V##1V##2V{\llap{$\vcenter{\hbox
 {$\scriptstyle##1$}}$}\Big\downarrow\rlap{$\vcenter{\hbox{%
$\scriptstyle##2$}}$}&}%
\def\={&\hskip.5em\mathrel
 {\vbox{\hrule width\minCDaw@\vskip3\ex@\hrule width
 \minCDaw@}}\hskip.5em&}
\def\verteq{\Big\Vert&}
\def\novarr{&}
\def\noharr{&&}
\def\SE##1E##2E{\slantedarrow(0,18)(4,-3){##1}{##2}&}
\def\SW##1W##2W{\slantedarrow(24,18)(-4,-3){##1}{##2}&}
\def\NE##1E##2E{\slantedarrow(0,0)(4,3){##1}{##2}&}
\def\NW##1W##2W{\slantedarrow(24,0)(-4,3){##1}{##2}&}
\def\slantedarrow(##1)(##2)##3##4{%
\thinlines\unitlength1pt\lower 6.5pt\hbox{\begin{picture}(24,18)%
\put(##1){\vector(##2){24}}%
\put(0,8){$\scriptstyle##3$}%
\put(20,8){$\scriptstyle##4$}%
\end{picture}}}
\def\vspace##1{\noalign{\vskip##1\relax}}\relax\let\amp@rs@nd@&\iffalse}\fi
 \CD@true\vcenter\bgroup\relax\let\\=\cr\iffalse}\fi\tabskip\z@skip\baselineskip20\ex@
 \lineskip3\ex@\lineskiplimit3\ex@\halign\bgroup
 &\hfill$\m@th##$\hfill\cr}
\def\@endCDS{\cr\egroup\egroup}
\newdimen\TriCDarrw@
\newif\ifTriV@
\newenvironment{TriCDV}{\@TriCDV}{\@endTriCD}
\newenvironment{TriCDA}{\@TriCDA}{\@endTriCD}
\def\@TriCDV{\TriV@true\def\TriCDpos@{6}\@TriCD}
\def\@TriCDA{\TriV@false\def\TriCDpos@{10}\@TriCD}
\def\@TriCD#1#2#3#4#5#6{%
\setbox0\hbox{$\ifTriV@#6\else#1\fi$}
\TriCDarrw@=\wd0 \advance\TriCDarrw@ 24pt
\advance\TriCDarrw@ -1em
\def\SE##1E##2E{\slantedarrow(0,18)(2,-3){##1}{##2}&}
\def\SW##1W##2W{\slantedarrow(12,18)(-2,-3){##1}{##2}&}
\def\NE##1E##2E{\slantedarrow(0,0)(2,3){##1}{##2}&}
\def\NW##1W##2W{\slantedarrow(12,0)(-2,3){##1}{##2}&}
\def\slantedarrow(##1)(##2)##3##4{\thinlines\unitlength1pt
\lower 6.5pt\hbox{\begin{picture}(12,18)%
\put(##1){\vector(##2){12}}%
\put(-4,\TriCDpos@){$\scriptstyle##3$}%
\put(12,\TriCDpos@){$\scriptstyle##4$}%
\end{picture}}}
\def\={\mathrel {\vbox{\hrule
   width\TriCDarrw@\vskip3\ex@\hrule width
   \TriCDarrw@}}}
\def\>##1>>{\setbox\z@\hbox{$\scriptstyle
 \;{##1}\;\;$}\global\bigaw@\TriCDarrw@
 \ifdim\wd\z@>\bigaw@\global\bigaw@\wd\z@\fi
 \hskip.5em
 \mathrel{\mathop{\hbox to \TriCDarrw@
{\rightarrowfill}}\limits^{##1}}
 \hskip.5em}
\def\<##1<<{\setbox\z@\hbox{$\scriptstyle
 \;{##1}\;\;$}\global\bigaw@\TriCDarrw@
 \ifdim\wd\z@>\bigaw@\global\bigaw@\wd\z@\fi
 \mathrel{\mathop{\hbox to\bigaw@{\leftarrowfill}}\limits^{##1}}
 }
 \CD@true\vcenter\bgroup\relax\let\\=\cr\iffalse}\fi
 \tabskip\z@skip\baselineskip20\ex@
 \lineskip3\ex@\lineskiplimit3\ex@
 \ifTriV@
 \halign\bgroup
 &\hfill$\m@th##$\hfill\cr
#1&\multispan3\hfill$#2$\hfill&#3\\
&#4&#5\\
&&#6\cr\egroup%
\else
 \halign\bgroup
 &\hfill$\m@th##$\hfill\cr
&&#1\\%
&#2&#3\\
#4&\multispan3\hfill$#5$\hfill&#6\cr\egroup
\fi}
\def\@endTriCD{\egroup}
\newcommand{\sA}{{\mathcal A}}
\newcommand{\sB}{{\mathcal B}}
\newcommand{\sC}{{\mathcal C}}
\newcommand{\sD}{{\mathcal D}}
\newcommand{\sE}{{\mathcal E}}
\newcommand{\sF}{{\mathcal F}}
\newcommand{\sG}{{\mathcal G}}
\newcommand{\sH}{{\mathcal H}}
\newcommand{\sI}{{\mathcal I}}
\newcommand{\sJ}{{\mathcal J}}
\newcommand{\sK}{{\mathcal K}}
\newcommand{\sL}{{\mathcal L}}
\newcommand{\sM}{{\mathcal M}}
\newcommand{\sN}{{\mathcal N}}
\newcommand{\sO}{{\mathcal O}}
\newcommand{\sP}{{\mathcal P}}
\newcommand{\sQ}{{\mathcal Q}}
\newcommand{\sR}{{\mathcal R}}
\newcommand{\sS}{{\mathcal S}}
\newcommand{\sT}{{\mathcal T}}
\newcommand{\sU}{{\mathcal U}}
\newcommand{\sV}{{\mathcal V}}
\newcommand{\sW}{{\mathcal W}}
\newcommand{\sX}{{\mathcal X}}
\newcommand{\sY}{{\mathcal Y}}
\newcommand{\sZ}{{\mathcal Z}}
\newcommand{\A}{{\mathbb A}}
\newcommand{\B}{{\mathbb B}}
\newcommand{\C}{{\mathbb C}}
\newcommand{\D}{{\mathbb D}}
\newcommand{\E}{{\mathbb E}}
\newcommand{\F}{{\mathbb F}}
\newcommand{\G}{{\mathbb G}}
\newcommand{\bH}{{\mathbb H}}
\newcommand{\I}{{\mathbb I}}
\newcommand{\J}{{\mathbb J}}
\newcommand{\bL}{{\mathbb L}}
\newcommand{\M}{{\mathbb M}}
\newcommand{\N}{{\mathbb N}}
\newcommand{\bP}{{\mathbb P}}
\newcommand{\Q}{{\mathbb Q}}
\newcommand{\R}{{\mathbb R}}
\newcommand{\T}{{\mathbb T}}
\newcommand{\U}{{\mathbb U}}
\newcommand{\V}{{\mathbb V}}
\newcommand{\W}{{\mathbb W}}
\newcommand{\X}{{\mathbb X}}
\newcommand{\Y}{{\mathbb Y}}
\newcommand{\Z}{{\mathbb Z}}
\newcommand{\rk}{{\rm rank}}
\title[Subvarieties of moduli stacks]{Discreteness of minimal models of Kodaira dimension zero
and subvarieties of moduli stacks}
\author[Eckart Viehweg]{Eckart Viehweg}
\address{Universit\"at Essen, FB6 Mathematik, 45117 Essen, Germany}
\email{ viehweg@uni-essen.de}
\thanks{This work has been supported by the ``DFG-Schwerpunktprogramm
Globale Methoden in der Komplexen Geometrie''. The second named author
is supported by a grant from the Research
Grants Council of the Hong Kong
Special Administrative Region, China
(Project No. CUHK 4239/01P)}
\author[Kang Zuo]{Kang Zuo}
\address{The Chinese University of Hong Kong,
Department of Mathematics,
Shatin, Hong Kong}
\email{kzuo@math.cuhk.edu.hk}
\maketitle
\tableofcontents
\section{Introduction}
Let $f:V\to U$ be a smooth projective morphism with connected
fibres over a complex quasi-projective manifold $U$.
\begin{definition}\label{variation} \
\begin{enumerate}
\item[i.] ${\rm Var}(f)$ is the smallest integer $\eta$ for which there
exists a finitely generated subfield $K$ of $\overline{\C(U)}$ of transcendence degree
$\eta$ over $\C$, a variety $F'$ defined over $K$, and a
birational equivalence
$$
V\times_U {\rm Spec}(\overline{\C(U)}) \sim F'\times_{{\rm Spec}(K)}
{\rm Spec}(\overline{\C(U)}).
$$
\item[ii.] $f:V\to Y$ is birationally isotrivial if ${\rm Var}(f)= 0$, hence if
there exists some generically finite covering $U'\to U$, a projective manifold $F'$,
and a birational map
$$
V\times_UU' \sim U'\times F'.
$$
\item[iii.]
$f:V\to U$ is (biregulary) isotrivial if there exists a generically finite
covering $U'\to U$, a projective manifold $F$ and an isomorphism
$$
V\times_UU' \simeq U'\times F.
$$
\end{enumerate}
\end{definition}
So the variation of a morphisms counts the number of parameters
controlling the birational structure of the fibres of $F$.

Maehara has shown in \cite{Mae} that under the assumption that
$\omega_F$ is semi-ample and big for a general fibre $F$ of $f$, a family is
birationally isotrivial, if and only if it is biregulary isotrivial.
In different terms, for families of minimal models of complex manifolds
of general type, ${\rm Var}(f)$ measures the number of directions
where the structure of $F$ varies. As shortly discussed in \ref{sectcan}
his result today follow immediately from the existence of the moduli scheme
$M_h$ of canonically polarized manifolds, and from the description of an
ample invertible sheaf on $M_h$.

We will slightly extend the methods used to prove \cite{Vie},
Theorem 6.24, to show that for families with
$\omega^\delta_F=\sO_F$, for some $\delta >0$, the same holds
true. We will show that for a given projective manifold $F'$ the
set of minimal models is discrete, hence that there are no
non-trivial families of minimal models.

Let us fix some polarization $\sL$ of $f:V\to U$, with Hilbert polynomial $h$.
If $\omega_{V/U}$ is $f$-ample we will choose $\sL=\omega_{V/U}^\rho$,
for some $\rho >0$. By \cite{Vie} there exists a quasi-projective moduli schemes $M_h$,
parameterizing polarized manifolds $(F,\sL)$ with $\omega_F$ semiample
and with $h(\nu)=\chi(\sL^\nu)$.
The family $f:V\to U$ together with $\sL$ induces a map
$$\varphi:U \to M_h$$
Since we require $\varphi:U\to M_h$
to be induced by a family it factors through the moduli stack
$\sM_h$.
\begin{theorem}\label{stthm}
Let $f:V\to U$ be a family of polarized manifolds.
Assume that $\omega_F^\delta=\sO_F$, for some $\delta >0$
(or that all fibres $F$ of $f$ are canonically polarized). Let
$\varphi:U\to M_h$ be the induced morphism to the moduli
scheme. Then ${\rm Var}(f)=\dim(\varphi(U))$.
\end{theorem}
If in Theorem \ref{stthm} the morphism $f:V\to U$ is birationally isotrivial,
$\varphi(U)$ must be zero dimensional, hence $f$ is biregulary isotrivial.

As Y. Kawamata told us, Theorem \ref{stthm} remains true for
families of polarized manifolds with $\omega_{X/Y}$ $f$-semiample.
Here however one has to replace the moduli scheme $M_h$ by the
moduli scheme $P_h$ of polarized manifolds, up to numerical
equivalence (see \cite{Vie}). So given a family $f:V\to U$ of
polarized manifolds with $\omega_{V/U}$ $f$-semiample, let $\psi:U
\to P_h$ be the morphism to the moduli scheme of polarized
manifolds up to numerical equivalence. Then ${\rm
Var}(f)=\dim(\psi(U))$.

For some applications to subvarieties of moduli stacks of polarized
$n$-folds of Kodaira dimension $0<\kappa<n$, one still has to
understand the structure of the ample sheaves on $P_h$, less
transparent than the ones on $M_h$. Nevertheless, we hope that most of the results
stated in the second half of this note remain true for all Kodaira dimensions,
with $M_h$ replaced by $P_h$.

Theorem \ref{stthm}, or more generally the equivalence of
biregular and birational isotriviality allows to extends some of the results
obtained in \cite{VZ3} for canonically polarized manifolds to
families of manifolds $F$ with $\omega_F^\delta=\sO_F$ (see
Theorem \ref{stthm2}, i) and ii) and Section \ref{sectappl}).
This is done in the second half of this article, a continuation of \cite{VZ3}.
What methods are concerned, the reader familiar with \cite{VZ3}
will find nothing new. In fact we just
sketch the changes needed to extend some of the results to this
case.

In the final
Section \ref{rigidity} we will state a criterion for the rigidity of
non-isotrivial families over curves, and its translation to curves
in the moduli stack of minimal polarized manifolds of Kodaira
dimension zero, or of canonically polarized manifolds. This
criterion is implicitly used in \cite{VZ3}, Proof of 6.4 and 6.5,
but it was not explicitly stated there.

A slightly weaker statement (\ref{proprigidsa}) extends to all families
with $\omega_F$ semiample. A similar
criterion has been shown by S. Kov\'acs and, for families of Calabi-Yau
manifolds by  K. Liu, A. Todorov, S.-T. Yau and the second named author
in \cite{LTYZ}. As a corollary one obtains (see \ref{corrigid2}):

\begin{corollary}\label{corrigid3} \
Let $M_h$ be either the moduli scheme of canonically polarized
manifolds or the moduli scheme of polarized manifolds $F$ with
$\omega_F^\delta=\sO$ for some $\delta >0$.
There are only finitely many morphisms $\varphi:U\to M_h$
which are induced by a smooth family $f:V\to U$ with:\\[.2cm]
For a general fibre $F$ of $f$ the $n$-th wedge product
$$0\neq \wedge^{n} \xi \in H^{n}(F,\omega^{-1}_F),$$
where $\xi\in H^1(F,T_F)$ denotes the
Kodaira Spencer class corresponding to the deformation $f:V\to U$ of $F$.
\end{corollary}

The first half of this article presents a proof of Theorem \ref{stthm},
hopefully of interest independently of the applications to subvarieties
of moduli stacks. In the first section, we will show, that the proof of
Theorem \ref{stthm} can be reduced to families over a curve. Next we recall
and strengthen a Positivity Theorem from \cite{Vie}. It allows to reinterpret
Maehara's Result in section \ref{sectcan}. The case of minimal models of
Kodaira dimension zero is handled in Section~\ref{sectzero}.

\section{Reduction to families over curves}\label{sectred}
In Definition \ref{variation}, i), we may choose a finitely generated subfield
$L$ of $\overline{\C(U)}$ which contains $\C(U)$ and $K$. Let $U'$ be the normalization
of $U$ in $L$, and let $T$ be a smooth quasi-projective variety with function field
$\C(T)=K$. Replacing $T$ by some open subscheme, we may assume that there exists
a smooth projective morphism $g:Z\to T$ with general fibre $F'$, and replacing
$U'$ by some open subscheme one finds morphisms $\tau:U'\to U$ and
$\pi:U'\to T$ fitting into a diagram
$$
\begin{CD}
V \<<< V' \>\sim >> Z' \>>> Z\\
\V f VV \V f' VV \V g' VV \V g VV\\
U \< \tau << U' \> = >> U' \> \pi >> T,
\end{CD}
$$
where $V'\to Z'$ is a birational equivalence, and where the right and left hand
squares are fibre products. For a point $t \in T$ in general position, one has
$$
\dim(\pi^{-1}(t))= \dim(U')-\dim(T) = \dim(U) - {\rm Var}(f).
$$
If under the assumptions made in Theorem \ref{stthm} ${\rm
Var}(f) < \dim(\varphi(U))$, for a point $\eta \in \varphi(U)$ in
general position,
$$
\dim(\tau^{-1}\varphi^{-1}(\eta))= \dim(U')-\dim(\varphi(U))
< \dim(\pi^{-1}(t)),
$$
hence there exists a curve $C$ in $\pi^{-1}(t)$ with
$\tau\circ\varphi|_C$ finite. In order to prove Theorem
\ref{stthm} one just has to show that such a curve can not
exist. Theorem \ref{stthm} follows from
\begin{proposition}\label{mainprop} Let $U$ be a non-singular irreducible
curve and let $F'$ be a projective manifold. Let $f:V\to U$ be a family
of polarized manifolds. Assume that there
exists a birational equivalence $V \sim U\times F'$ over $U$.
If either the fibres $F$ of $f$ are canonically polarized, or if
$\omega_F^\delta=\sO_F$, for some $\delta >0$, then the induced morphism
$\varphi:U\to M_h$ is constant.

In particular, $f:V\to U$ is biregulary isotrivial.
\end{proposition}
In order to prove Proposition \ref{mainprop} we may replace $U$ by a finite covering.
Doing so one can assume that $f:V\to U$ extends to a semistable morphism
$f:X\to Y$ of projective manifolds, hence that $\Delta=f^{-1}(S)$ is a reduced
normal crossing divisor, for $S=Y\setminus U$.

Moreover, we may replace the given polarization by some power. In fact, the corresponding
map of the moduli schemes $M_h$ is a finite map. This allows to assume that
for all fibres $F$ of $V\to U$ the polarization $\sL$ is very ample, and without higher
cohomology.

\section{Positivity of direct image sheaves}\label{sectpos}

Recall that a locally free sheaf $\sE$ on a projective non-singular
curve $Y$ is numerically effective (nef),
if for all finite morphisms $\tau:Z\to Y$ and for all invertible quotients
$\sL$ of $\tau^*(\sE)$ the degree $\deg(\sL)\geq 0$.

Fujita's positivity theorem (today an easy corollary of
Koll\'ar's vanishing theorem) says that $f_* \omega_{X/Y}$ is
nef. By \cite{VZ1}, 2.3, one obtains as a direct consequence.
\begin{lemma} \label{fujita} Let $f:X\to Y$ be a morphism from a
normal projective variety $X$ to a curve $Y$, with connected
fibres. Assume that $X$ has at most rational double points as singularities.
Let $\sN$ be an invertible sheaf on $X$ and $\Gamma$ an effective divisor. Assume
that for some $N > 0$ there exists a nef locally free sheaf $\sE$ on $Y$ and a
surjection $$f^* \sE \>>> \sN^N (- \Gamma).$$
Then
$$f_* \left(\sN \otimes \omega_{X/Y} \left\{ - \frac{\Gamma}{N}
\right\} \right)$$
is nef.
\end{lemma}
Here $\omega_{X/Y} \left\{ - \frac{\Gamma}{N} \right\}$ denotes
the (algebraic) multiplier sheaf (see for example \cite{EV}, 7.4,
or \cite{Vie}, section 5.3). If $\tau: X' \to X$ is any blowing up
with $\Gamma' = \tau^* \Gamma$ a normal crossing divisor, then $$
\omega_{X/Y} \left\{ - \frac{\Gamma}{N} \right\} = \tau_*
\left(\omega_{X'/Y} \left( - \left[ \frac{\Gamma'}{N} \right]
\right)\right). $$ As in \cite{EV}, \S \ 7 and \cite{Vie}, section
5.3, we are mainly interested in the case where the multiplier
sheaf on a general fibre $F$ is isomorphic to $\omega_F$. The
corresponding threshold is defined for any effective divisor $\Pi$
or any invertible sheaf $\sL$ on $F$ with $H^0 (F, \sL) \neq 0$.
\begin{gather*}
e(\Pi) = {\rm Min} \left\{ N \in \N - \{ 0 \} ; \ \omega_F \left\{-
\frac{\Pi}{N} \right\} = \omega_F \right\} \ \ \mbox{ \ \ \
\ \ and}\\
e (\sL) = {\rm Max} \left\{ e (\Pi); \ \Pi \ \mbox{the zero set of} \
\sigma \in H^0 (F, \sL) - \{ 0 \} \right\}.
\end{gather*}

For smooth morphisms $f:X\to Y$ and for an $f$-ample sheaf $\sL$ on $X$
we obtained in \cite{Vie}, 6.24 and 7.20, strong positivity theorems.
Their proof, in case $Y$ is a curve, can easily be extended to
semistable morphisms $f:X\to Y$.
\begin{theorem}\label{positiv}
Let $Y$ be a curve, let $f:X\to Y$ be a semistable morphism
between projective manifolds with connected fibres, and let $\sM$ be an
invertible sheaf on $X$. Let $U\subset Y$ be an open dense
subscheme with $V=f^{-1}(U) \to U$ smooth.
Assume that for all fibres $F$ of $V\to U$ the canonical sheaf $\omega_F$
is semiample, that $\sM|_F$ is very ample and without higher cohomology.
Then for
\begin{gather*}
e \geq c_1(\sM|_F)^{\dim(F)}+2, \ \ \ \ \
r=\rk (f_*\sM)\\
\mbox{ and \ \ \ }r(\nu)=\rk (f_*(\sM^\nu\otimes \omega^{e\cdot\nu}_{X/Y})):
\end{gather*}
\begin{enumerate}
\item[a.] For all $\nu >0$ $$\big( \bigotimes^r f_*(\sM^\nu\otimes
\omega_{X/Y}^{e\cdot\nu})\big)
\otimes\det(f_*\sM)^{-\nu}$$
is nef.
\item[b.] If the invertible sheaf
$\det(f_*(\sM^\nu\otimes\omega^{e\cdot\nu}_{X/Y}))\otimes\det(f_*\sM)^{-\nu
r(\nu)}$ is ample for some $\nu>0$,
$$\big( \bigotimes^r f_*(\sM\otimes \omega_{X/Y}^{e})\big)
\otimes\det(f_*\sM)^{-1}$$ is ample.
\item[c.] If for all $\nu > 0$ the degree of
$$\det(f_*(\sM^\nu\otimes\omega^{e\cdot\nu}_{X/Y}))\otimes\det(f_*\sM)^{-\nu
r(\nu)}$$ is zero, then $f:V\to U$ is biregulary isotrivial as a family of
polarized manifolds, i.e.
there exists some finite covering $U'\to U$, a projective
manifold $F'$, invertible sheaves $\sL'$ on $F'$ and $\sB$ on $U'$,
and an isomorphism
$$
\pi:V'=X\times_YU' \to F'\times U'
$$
with
$$
pr_1^*\sM = \pi^*( pr_1^*\sL\otimes pr_2^*\sB).
$$
\end{enumerate}
\end{theorem}
\begin{proof} As indicated already, the proof of
parts a) and b) will follow the arguments used in \cite{Vie}, 194--196, to prove
6.20. We just have to take care, that for a semistable family over a curve
the sheaves are at most getting larger. So we repeat the arguments.

For a) let us fix some $\nu >0$. For b) we assume that
$$
\det(f_*(\sM^\nu\otimes\omega^{e\cdot\nu}_{X/Y}))
\otimes\det(f_*\sM)^{-\nu \cdot r(\nu)}
$$
is ample.

The semicontinuity of the threshold, shown in \cite{Vie}, 5.17 for example,
allows to find some $\gamma \geq e\cdot\nu$ with
\begin{equation}\label{gamma}
e (\sM|_F ^{\nu \cdot e} \otimes
\omega^{e \cdot \nu \cdot (e-1)}_{F} ) \leq \gamma
\end{equation}
for all fibres $F$ of $V \to U$.

For b) we will show that
\begin{gather*}
S^{\gamma} \big(\big( \bigotimes^{r\cdot r(\nu)} f_*(\sM\otimes \omega_{X/Y}^{e})\big)
\otimes\det(f_*\sM)^{-{r(\nu)}}\big)\otimes\\
\det(f_*(\sM^\nu\otimes\omega^{e\cdot\nu}_{X/Y}))^{-r\cdot(e-1)}
\otimes\det(f_*\sM)^{\nu\cdot
r\cdot(e-1) \cdot r(\nu)}
\end{gather*}
is nef. Hence for both, a) or b), it is sufficient to prove
the corresponding statements for the pullback of the sheaves to any finite
covering $Y'$ of $Y$. Since we assumed $X\to Y$ to be semistable, the fibre
product $X'=X\times_YY'$ is a normal variety with at most rational double
points. Flat base change allows to replace $Y$ by such a covering and $(f:X\to
Y,\sM)$ by a desingularization of the pullback family.

Doing so, we may assume that $\det(f_*\sM)$ is the $r$-th power of an invertible
sheaf, and since all the sheaves occurring in a), b) or c) are compatible
with
changing the polarization by the pullback of an invertible sheaf on $Y$, we can
as well assume that $\det(f_*\sM)=\sO_Y$. Under this additional assumption
we have to verify in a) that
$$
f_*(\sM^\nu\otimes \omega_{X/Y}^{e\cdot\nu})
$$
is nef. For part b) we may assume in addition
that
$$
\det(f_* (\sM^{\nu} \otimes \omega^{e \cdot \nu}_{X/Y} ))^{r\cdot (e-1)}
= \sO_Y(\gamma \cdot H)
$$
for some effective divisor $H$ supported in $U$. We have to prove that
$$
\big( \bigotimes^{r\cdot r(\nu)} f_*(\sM\otimes \omega_{X/Y}^{e})\big)
\otimes \sO_Y(-H)
$$
is nef.\\

Let $f^s:X^s \to Y$ be the $s$-fold fibre product. $X^s$ is normal
with at most rational double points (see \cite{Mor}, page 291, for
example). Consider
$$
\sP = \bigotimes^{s}_{i=1} pr^{*}_{i} \sM .
$$
By flat base change one obtains
$$
f^{s}_{*}( \sP^{\alpha} \otimes \omega^{\beta}_{X^s /Y}) =
\bigotimes^s f_* (\sM^{\alpha} \otimes \omega^{\beta}_{X/Y})
$$
for all $\alpha , \beta$. The restriction of $\sP^{\nu} \otimes
\omega^{e\cdot\nu - \iota}_{X^s /Y} $ to ${f^s}^{-1}(U)=V^s$ is
$f^s$-ample for all $\iota \leq e\cdot\nu $. Let us write $\epsilon= e\cdot\nu$
or $\epsilon=e\cdot\nu -1$, where $\nu$ may be any positive integer.

If $\Gamma$ is the zero divisor
of a section of $\sP$, which does not contain any
fibre $F^s$ of $V^s\to U$ the compatibility of the threshold
with products and its semicontinuity imply (see \cite{Vie}, 5.14 and 5.21)
\begin{equation}\label{threshold}
e (\Gamma |_{F^s} ) \leq e (\sP |_{F^s} ) =
e (\sM |_{F} ) < e \mbox{ \ \ and \ \ }
e (\Gamma|_{V^s} ) < e.
\end{equation}
In fact, as shown in \cite{Vie}, 5.11), one has
$$
e(\sM|_F) \leq c_1(\sM|_F)^{\dim(F)} + 1.
$$
Moreover, by the choice of $\epsilon$
\begin{equation}\label{threshold2}
e(\nu\cdot\Gamma|_{F^s}) \leq \nu\cdot e(\Gamma|_{F^s}) \leq
\nu \cdot (c_1(\sM|_{F})^{\dim(F)} + 1) \leq \epsilon.
\end{equation}
Let $\sH$ be an ample invertible sheaf on $Y$.
\begin{claim}\label{claim}
Assume that for some $\rho \geq 0$, $N >0$, $M_0 >0$ and for all
multiples $M$ of $M_0$, the sheaf
$$
f_* ((\sM^{\nu} \otimes \omega^{\epsilon}_{X/Y} )^{M\cdot N} ) \otimes
\sH^{\rho \cdot \epsilon \cdot N \cdot M}
$$
is nef. Then
$$
f_* ((\sM^{\nu} \otimes \omega^{\epsilon}_{X/Y} )^N ) \otimes \sH^{\rho
\cdot (\epsilon \cdot N -1)}
$$
is nef.
\end{claim}

\begin{proof}
Let us choose $s=r$. The determinant gives an inclusion
$$
{\rm det} (f_* \sM) = \sO_Y \>>> f^{r}_{*} \sP
= \bigotimes^r f_* \sM,
$$
which splits locally. Hence the zero divisor $\Gamma$ of the
induced section of $\sP$ does not contain any fibre of
$V^r\to U$. For
$$
\sN = \sP^{\nu \cdot N} \otimes \omega^{\epsilon\cdot N -1}_{X^r /Y}
\otimes f^{r *} \sH^{\rho \cdot (\epsilon \cdot N -1) \cdot r}
$$
one obtains that the restriction of
$$
\sN^{\epsilon} (-\nu \cdot \Gamma ) = (\sP^{\nu} \otimes
\omega^{\epsilon}_{X^r /Y} \otimes f^{r*} \sH^{\rho \cdot \epsilon \cdot r}
)^{(\epsilon \cdot N -1)}
$$
to $V^s$ is $f^r$-ample.
If $M'$ is a positive integer, divisible by $M_0 \cdot N$ the
sheaf
$$
f^{r}_{*} (\sN^{\epsilon} (-\nu \cdot \Gamma )^{M'}) =
\bigotimes^r (f_* (\sM^{\nu} \otimes \omega^{\epsilon}_{X/Y} )^{(\epsilon
\cdot N -1) \cdot M'} \otimes \sH^{\rho \cdot \epsilon
\cdot r (\epsilon \cdot N -1) \cdot M'} )
$$
is nef. Choose $M'$ such that
$$
f^*f_*((\sM^\nu\otimes\omega^\epsilon_{X/Y})^{(\epsilon
\cdot N -1) \cdot M'}) \>>> (\sM^\nu\otimes\omega^\epsilon_{X/Y})^{(\epsilon
\cdot N -1) \cdot M'}
$$
is surjective over $U$.

\ref{fujita} implies that the subsheaf $f^{r}_{*} (\sN \otimes
\omega_{X^r /Y} \left\{ -\frac{\nu\cdot\Gamma}{\epsilon}\right\})$
of
$$
f^{r}_{*} (\sN \otimes \omega_{X^r /Y} ) = \bigotimes^r (f_* (
\sM^{\nu \cdot N} \otimes \omega^{\epsilon \cdot N}_{X/Y} ) \otimes
\sH^{\rho \cdot (\epsilon \cdot N -1)} )
$$
is nef. On the other hand, (\ref{threshold}) and (\ref{threshold2}) imply that
both sheaves coincide on $U$.
\end{proof}

Choose some $N_0>0$ such that for all multiples $N$ of $N_0$ and
for all $M>0$ the multiplication maps
$$
m: S^M (f_* (\sM^{\nu \cdot N} \otimes \omega^{\epsilon \cdot N}_{X/Y}
)) \>>> f_* (\sM^{\nu \cdot N \cdot M} \otimes \omega^{\epsilon \cdot N
\cdot M}_{X/Y} )
$$
are surjective over $U$. Define
$$
\rho = {\rm Min } \{ \mu >0 ; \ f_* (\sM^{\nu \cdot N} \otimes
\omega^{\epsilon \cdot N}_{X/Y} ) \otimes \sH^{\mu \cdot \epsilon \cdot N} \
\mbox{is nef} \} .
$$
The surjectivity of $m$ implies that
$$
f_* (\sM^{\nu \cdot N \cdot M} \otimes \omega^{\epsilon \cdot N \cdot
M}_{X/Y} ) \otimes \sH^{\rho \cdot \epsilon \cdot N \cdot M}
$$
is nef for all $M >0$. By \ref{claim}
$$
f_* (\sM^{\nu \cdot N} \otimes \omega^{\epsilon \cdot N}_{X/Y} )
\otimes \sH^{\rho \cdot (\epsilon \cdot N -1)} .
$$
is nef, hence by the choice of $\rho$
$$
(\rho -1)\cdot \epsilon \cdot N < \rho \cdot (\epsilon \cdot N -1)
$$
or equivalently $\rho < \epsilon \cdot N$. Then
$$
f_* (\sM^N \otimes \omega^{\epsilon \cdot N}_{X/Y} ) \otimes \sH^{\epsilon^2
\cdot N^2}
$$
is nef. This remains true if one replaces $Y$ by any finite covering, and
by \cite{VZ1}, 2.2, one obtains
that $f_* (\sM^N \otimes \omega^{\epsilon \cdot N}_{X/Y} )$ is nef.
Applying \ref{claim} a second time,
for the numbers $(N',N_0)$ instead of $(N,M_0)$ and for $\rho =0$,
one finds
\begin{claim}\label{nef2} For $\nu >0$ and $\epsilon=e\cdot\nu$ or
$\epsilon=e\cdot\nu -1$ and for all $N'>0$ the sheaf
$$
f_* ((\sM^{\nu} \otimes \omega^{\epsilon}_{X/Y} )^{N'})
$$
is nef.
\end{claim}
In particular, choosing $N'=1$ and $\epsilon=\nu e$ one obtains a).\\

For b) we consider the $s$-fold product $f^s : X^s \to Y$
for $s = r \cdot r (\nu).$
One has natural inclusions, splitting locally,
$$
\sO_Y = {\rm det} (f_* \sM)^{r (\nu)} \>>> f^{s}_{*} \sP
= \bigotimes^s f_*\sM
$$
and
$$
\det (f_* (\sM^{\nu} \otimes \omega^{e \cdot
\nu}_{X/Y} ))^{r} \>>> f^{s}_{*} (\sP^{\nu} \otimes
\omega^{e \cdot \nu}_{X^s /Y} ) = \bigotimes^s f_*
(\sM^{\nu} \otimes \omega^{e \cdot \nu}_{X/Y} ).
$$
If $\Delta_1$ and $\Delta_2$ denote the corresponding
zero-divisors on $X^s$ then $\Delta_1 + \Delta_2$ does not
contain any fibre of $V^s \to U$. Then
$$
\sP^{e\cdot\nu}\otimes\omega_{X^s/Y}^{e\cdot\nu\cdot(e-1)}
= {f^s}^*\det (f_* (\sM^{\nu} \otimes \omega^{e \cdot
\nu}_{X/Y} ))^{r\cdot(e-1)}\otimes \sO_{X^s}((e -1) \cdot \Delta_2 + \nu \cdot \Delta_1 ),
$$
and
$$
\sP^{\gamma} \otimes \omega^{\gamma \cdot (e -1)}_{X^s/Y}
= (\sP\otimes \omega_{X^s/Y}^{e-1})^{\gamma-\nu e}\otimes \sO_X( \gamma \cdot
{f^s}^* H +(e -1) \cdot \Delta_2 + \nu \cdot \Delta_1)).
$$
By \ref{nef2} the sheaf
$$
{f^s}_* ((\sP\otimes \omega_{X^s/Y}^{e-1})^{\gamma-\nu \cdot e})^M
= \bigotimes^s f_* ((\sP\otimes \omega_{X/Y}^{e-1})^{\gamma-\nu \cdot e})^M
$$
is nef for all $M>0$. \ref{fujita} implies that
$$
\sP \otimes \omega^{e-1}_{X^s/Y} \otimes \omega_{X^s/Y}\left\{-
\frac{\gamma \cdot{f^s}^*H +(e -1) \cdot \Delta_2 + \nu \cdot \Delta_1}{\gamma}
\right\}
$$
is nef. By (\ref{threshold}) and (\ref{gamma})
\begin{gather*}
e(((e -1) \cdot \Delta_2 + \nu \cdot \Delta_1) |_{F^{s}} ) \leq e (\sP|_{F^{s}}^{\nu
\cdot e } \otimes \omega^{e \cdot \nu \cdot
(e -1)}_{F^{s}} ) =\\ e (\sM|_{F}^{\nu \cdot e} \otimes
\omega^{e \cdot \nu \cdot (e-1)}_{F} )\leq \gamma
\end{gather*}
for all fibres $F$ of $V \to U$. Hence the cokernel of
$$
\omega_{X^s/Y}\left\{
\frac{-\gamma \cdot {f^s}^* H +(e -1) \cdot \Delta_2 + \nu \cdot \Delta_1}{\gamma}
\right\} \to \omega_{X^s/Y}(-{f^s}^* H)
$$
lies in $X^s\setminus V^s$, and thereby
$$
{f^s}_* (\sP \otimes \omega^{e}_{X^s/Y}) \otimes \sO_Y(-H)=
\big( \bigotimes^{r\cdot r(\nu)} f_*(\sM\otimes \omega_{X/Y}^{e})\big)
\otimes \sO_Y(-H)
$$
is nef.\\

Part c) follows from part a) and Kollar's ampleness criterion (see \cite{Vie},
4.34). Again we may assume that
$\det(f_*\sM)=\sO_Y$. By part a) for all $\eta >0$ the sheaf $\sE=f_*(\sM^\eta \otimes
\omega_{X/Y}^{e\cdot\eta})$ is nef. Choose $\nu>0$ such that the multiplication map
$$
\mu: S^\nu(\sE^\eta) \to f_*(\sM^{\nu\cdot\eta}\otimes \omega_{X/Y}^{e\cdot\nu \cdot \eta})
$$
is surjective over $U$. By \cite{Vie}, 4.34, $\det({\rm Im}(\mu))$ is
ample, if the kernel $\sK$ of the multiplication map is of maximal
variation. Let us recall the definition. For a point $y\in U$ choose a local
trivialization of $\sE$. Then $\sK_y=\sK\otimes \C(y)$ as a subvectorspace of
$S^\nu(\C^{r(\eta)})$, defines a point $[\sK_y]$ in the Grassmann variety
$$
\G r ={\rm Grass}(r(\nu\cdot\eta), S^\nu(\C^{r(\eta)})).
$$
The group $G={\rm Sl}(r(\eta),\C)$ acts on $S^\nu(\C^{r(\eta)})$, hence on $\G r$.
Let $G_y$ denote the orbit of $[\sK_y]$. The kernel has maximal variation, if
$$\{ z \in Y; \ G_{z}=G_y\}$$ is finite, as well as the stabilizer of $[\sK_y]$.

The second condition holds true for $\eta$ and $\nu$ sufficiently large.
In fact, $\sK_y$ determines the fibre $f^{-1}(y)$ as a subvariety of
$$\bP(H^0(f^{-1}(y),(\sM^\eta
\otimes \omega_{X/Y}^{e\cdot\eta})|_{f^{-1}(y)})),
$$
and $[\sK_y]$ is nothing but the point of the Hilbert scheme $\rm Hilb$ parameterizing
subvarieties of this projective space. By \cite{Vie}, 7.2, the stabilizer
of such a point is finite.

The assumption in c) implies that $\sK$ is not of maximal variation.
Hence for all points $z$ in a neighborhood
$U_y$ of $y$, the orbits $G_z$ coincide. In different terms the images of $z\in U_y$ in
$\rm Hilb$ all belong to the same $G$-orbit. Since $M_h$ is a quotient of a subscheme of $\rm Hilb$
by the $G$-action, the morphism $\varphi:U \to M_h$ is constant, as claimed in c).
\end{proof}

\section{Families of canonically polarized manifolds}\label{sectcan}
For families with $\omega_F$ big
and semi-ample the equivalence of birational and biregular isotriviality has
been shown by Maehara in \cite{Mae}. For families of canonically polarized
manifolds, one just has to use, that the fibres are their own canonical model.

Or, to formulate the proof parallel to the one given below
in the Kodaira dimension zero case, one could argue in the following way.
Assume that $U$ is a curve, and choose a semistable compactification $f:X\to Y$ of
$V\to U$. By assumption $X$ is birational to the trivial family $F'\times Y$
over $Y$, hence $f_*\omega_{X/Y}^\nu$ is a direct sum of copies of $\sO_Y$, for
all $\nu$. Obviously, if $\sM$ is some power of $\omega_{X/Y}$ this implies that
$$
\det(f_*(\sM^\nu\otimes\omega^{e\cdot\nu}_{X/Y}))\otimes\det(f_*\sM)^{-\nu
r(\nu)}=\sO_Y,
$$
and by \ref{positiv}, c), one finds $V\to U$ to be biregulary isotrivial.

\section{Families of manifolds of Kodaira dimension zero}\label{sectzero}
Let $U$ be a curve and $f:V\to U$ be a family of polarized
manifolds $F$ with $\omega_F^\delta=\sO_F$. In order to prove
\ref{mainprop} we may replace $U$ by some finite cover, and we may
choose a compactification $f:X\to Y$ satisfying the assumptions
made in Theorem \ref{positiv}.

By assumption, $\lambda=f_*\omega_{X/Y}^{\delta}$ is an invertible sheaf,
and the natural map $f^*\lambda\to \omega_{X/Y}^\delta$ is an isomorphism
over $V$. Let $E$ be the zero divisor of this map, i.e.
$$
\omega_{X/Y}^\delta = f^*\lambda\otimes\sO_X(E).
$$
$E$ is supported in $\Delta=X\setminus V$ and, since the fibres of $f$ are reduced
divisors, $E$ can not contain a whole fibre. Hence for all $\mu >0$
$f_*\sO_X(\mu E)=\sO_Y$ and $f_*\omega_{X/Y}^{\mu\cdot\delta}=\lambda^\mu$.

Let $\sM'$ be any polarization which is very ample and without higher cohomology
on the fibres of $V\to U$. The sheaf
$$
f_*(\sM'\otimes\sO_X(*E))=f_*(\sM'\otimes(\lim_{\mu>0}\sO(\mu E)))
$$
is coherent, hence locally free.

In fact, locally \'etale or locally analytic
we can choose over a neighborhood $\sU$ of $s\in S=Y\setminus U$ a section
$\sigma:\sU \to X$ with image $C$, not meeting the support of $E$.
If $I$ denotes the ideal sheaf of $C$, for some $\rho \gg 0$
$$
f_*(\sM'|_{f^{-1}(\sU)}\otimes I^\rho)=0.
$$
Since direct images are torsion free, and since $E$ is supported in fibres,
$$
f_*((\sM'\otimes\sO_X(*E))|_{f^{-1}(\sU)}\otimes I^\rho)=0.
$$
Then $f_*(\sM'\otimes\sO_X(*E))|_\sU$ is a torsion free subsheaf of
$$
f_*((\sM'\otimes\sO_X(*E))|_{f^{-1}(\sU)}\otimes \sO_{f^{-1}(\sU)}/I^\rho)=
f_*(\sM'|_{f^{-1}(\sU)}\otimes \sO_{f^{-1}(\sU)}/I^\rho),
$$
hence coherent.

Let $\sM$ be the reflexive hull of the image of
$$
f^*f_*(\sM'\otimes\sO_X(*E))\>>> \sM'\otimes\sO_X(*E).
$$
$\sM$ is again coherent, and it must be contained in
$\sM'\otimes\sO_X(\alpha\cdot E)$ for some $\alpha$. Since it is
reflexive, it is an invertible sheaf. By construction
$\sM|_V\simeq \sM'|_V$ and
\begin{multline*}
f_*(\sM'\otimes\sO_X(*E))= f_*f^*f_*(\sM'\otimes\sO_X(*E)) \subset\\
f_*\sM \subset f_*(\sM\otimes\sO_X(*E)) \subset f_*(\sM'\otimes\sO_X(*E)),
\end{multline*}
hence all those sheaves coincide.
We found an invertible sheaf $\sM$ satisfying the assumptions made in
\ref{positiv} with the additional condition
\begin{equation*}
f_*(\sM\otimes \omega_{X/Y}^e)=f_*(\sM\otimes \sO_X(\frac{e}{\delta}\cdot E)
\otimes f^*\lambda^{\frac{e}{\delta}})=
f_*\sM\otimes \lambda^{\frac{e}{\delta}},
\end{equation*}
for all multiples $e$ of $\delta$. For those $e$
\begin{multline}\label{direct}
\big( \bigotimes^r f_*(\sM\otimes \omega_{X/Y}^{e})\big)
\otimes\det(f_*\sM)^{-1} =\\
\big( \bigotimes^r (f_*\sM)\otimes \lambda^{\frac{e}{\delta}}
\big) \otimes\det(f_*\sM)^{-1}.
\end{multline}
\begin{proof}[Proof of \ref{mainprop} for $\omega_F^\delta=\sO_F$] \ \\
By assumption $f:X\to Y$ is birational over $Y$ to the trivial family
$pr_2:F'\times Y \to Y$, hence
$$
\lambda=f_*\omega^\delta_{X/Y}=\sO_Y.
$$
So the sheaf in (\ref{direct}) is
$$
\big( \bigotimes^r (f_*\sM) \big) \otimes\det(f_*\sM)^{-1}.
$$
Since its determinant is of degree zero, it can not be ample.
(\ref{direct}) and \ref{positiv}, b), imply that for no $\nu >0$
the sheaf
$$
\det(f_*(\sM^\nu\otimes\omega^{e\cdot\nu}_{X/Y}))\otimes\det(f_*\sM)^{-\nu
\cdot r(\nu)}
$$
is ample. By \ref{positiv}, a), it is of non negative degree, and
\ref{positiv}, c), implies that the family $V\to U$ is biregulary isotrivial.
\end{proof}
\begin{remark} Assume that $\delta=1$, hence that
$\omega_{V/U}=f^*\lambda|_U$. The argument used in the proof of
\ref{mainprop} shows in this particular case that ``$f$~non-isotrivial''
implies that on the compactification $Y$ of $U$
$$
\deg(f_*\omega_{X/Y}) >0.
$$
Since the same holds true for all finite coverings of $U$, one
obtains that the fibres of the period map from $M_h$ to the period domain
classifying the corresponding variations of Hodge structures can
not contain a quasi-projective curve. Of course this is a well known
consequence of the local Torelli Theorem for manifolds with a trivial canonical
bundle.
\end{remark}
\section{Kodaira-Spencer maps}\label{kernels}
Recall first the following definition, replacing of nef and ample,
on projective manifolds $Y$ of higher dimension.
\begin{definition}\label{positivity}
Let $\sF$ be a torsion free coherent sheaf on a quasi-projective
normal variety $Y$ and let $\sH$ be an ample invertible sheaf.
\begin{enumerate}
\item[a)] $\sF$ is generically generated if the natural morphism
$$
H^0 (Y, \sF) \otimes \sO_Y \longrightarrow \sF
$$
is surjective over some open dense subset $U_0$ of $Y$. If one wants to
specify $U_0$ one says that $\sF$ is globally generated over $U_0$.
\item[b)] $\sF$ is weakly positive if there exists some
dense open subset $U_0$ of $Y$ with $\sF|_{U_0}$ locally free, and
if for all $\alpha > 0$ there exists some $\beta > 0$ such that
$$
S^{\alpha \cdot \beta} (\sF) \otimes \sH^{\beta}
$$
is globally generated over $U_0$. We will also say that $\sF$ is
weakly positive over $U_0$, in this case.
\item[c)] $\sF$ is big if there exists some open dense subset
$U_0$ in $Y$ and some $\mu > 0$ such that
$$
S^{\mu} (\sF) \otimes \sH^{-1}
$$
is weakly positive over $U_0$. Underlining the role of $U_0$ we will
also call $\sF$ ample with respect to $U_0$.
\end{enumerate}
\end{definition}
Here, as in \cite{Vie} and \cite{VZ3}, we use the following
convention: If $\sF$ is a coherent torsion free sheaf on a quasi-projective
normal variety $Y$, we consider the largest open subscheme $i: Y_1 \to
Y$ with $i^* \sF$ locally free. For
$$
\Phi = S^{\mu}, \ \ \
\Phi =\bigotimes^\mu \mbox{ \ \ \ or \ \ \ }\Phi = \det
$$
we define
$$
\Phi (\sF) = i_* \Phi (i^* \sF).
$$

Again, $f:V\to U$ denotes a smooth family of manifolds over a
quasi-projective manifold $U$, which is allowed to be of dimension
larger than one. We choose non-singular projective
compactifications $Y$ of $U$ and $X$ of $V$, such that both
$S=Y\setminus U$ and $\Delta=X\setminus V$ are normal crossing
divisors and such that $f$ extends to $f:X\to Y$. As usual $\eta$
will denote a closed point in sufficient general position on $U$
and $X_\eta$ the fibre of $f$ over $\eta$. We will write
$T^i_{X_\eta}$ (or $T^i_{X/Y}(-\log \Delta)$ $\ldots$) for the
$i$-th wedge product of $T_{X_\eta}$ (or of $T_{X/Y}(-\log
\Delta)=\Omega^1_{X/Y}(\log \Delta)^\vee$ $\ldots$).

Let $T_\eta$ denote the restriction $T_U\otimes \C$ of the tangent sheaf of $U$
to $\eta$. The Kodaira-Spencer map
$$T_\eta \>>> H^1(X_\eta,T_{X_\eta})$$
gives rise to
$$
\bigotimes^\nu T_\eta \>>> \bigotimes^\nu H^1(X_\eta,T_{X_\eta}) \>>>
H^\nu(X_\eta,T^\nu_{X_\eta}).
$$
The composite map factors through
$$
S^\nu(T_\eta)\>>> H^\nu(X_\eta,T^\nu_{X_\eta}).
$$
One defines
\begin{multline*}
\mu(f)={\rm Max}\{\nu\in\N-\{0\}; \ S^\nu(T_\eta)\>>> H^\nu(X_\eta,T^\nu_{X_\eta})
\mbox{ is non zero}\}.
\end{multline*}
Of course, $\mu(f) \leq n=\dim(X_\eta)$.
We do not know any criterion, implying that for
$f:V\to U$ one has $\mu(f)=\dim(V)-\dim(U)$.

For example, if $U$ is a curve and $f:V\to U$ a family of polarized manifolds,
restricting the tautological sequence to ${X_\eta}=f^{-1}(\eta)$
one obtains an extension
$$
0 \>>> T_{X_\eta} \>>> T_X|_{X_\eta} \>>> \mathcal O_{X_\eta} \>>> 0
$$
and the induced class $\xi_\eta\in H^1({X_\eta},T_{X_\eta})$. Then $\mu(f)=\mu$
if and only if for $\eta$ in general position, the wedge product
$\wedge^\mu \xi_\eta \in H^\mu({X_\eta},T_{X_\eta}^{\mu})$ is non-zero, whereas
$\wedge^{\mu+1} \xi_\eta \in H^{\mu+1}({X_\eta},T_{X_\eta}^{\mu+1})$ is zero.
\begin{problem}\label{problem2}
Are there properties of ${X_\eta}$ which imply that for all families
$V\to U$ over a curve $U$, with general fibre ${X_\eta}=f^{-1}(\eta)$ the class
$$\wedge^n \xi_\eta \in H^n({X_\eta},\omega_{X_\eta}^{-1})$$ is non-zero?
\end{problem}
Being optimistic, one could try in \ref{problem2} the condition
``$\Omega_{X_\eta}^1$ ample''.\\

A slight extension of the main result of \cite{VZ3} says:
\begin{theorem}\label{stthm2} Assume that for a general fibre $X_\eta$ of
$f:X \to Y$ either $\omega_{X_\eta}$ is ample, or
$\omega^\delta_{X_\eta}=\sO_{X_\eta}$, for some $\delta$.
\begin{enumerate}
\item[i.] Then for some $m >0$ the sheaf $S^m(\Omega_Y^1(\log S))$
contains an invertible subsheaf $\sM$ of Kodaira dimension ${\rm Var}(f)$.
\item[ii.] If ${\rm Var}(f)=\dim(U)$ the sheaf $S^{\mu(f)}(\Omega_Y^1(\log S))$
contains a big coherent subsheaf $\sP$.
\item[iii.] Let $Z$ be a submanifold of $Y$ such that $S_Z=S\cap Z$ remains a
normal crossing divisor, and such that
$W=X\times_YZ$ is non-singular. For the induced family
$h:W \to Z$ assume that $\mu(f)=\mu(h)$. Then,
if ${\rm Var}(f)=\dim(U)$, the restriction of the sheaf $\sP$ from part ii)
to $S^{\mu(f)}(\Omega^1_Z(\log S_Z))$ is non trivial.
\item[iv.] Assume in iii) that $h:W\to Z$ is a desingularization of
the pullback of a family
$h':W'\to Z'$ under $\pi:Z\to Z'$, with $Z'$ non-singular and with $h'$ smooth
over $Z'\setminus S_{Z'}$ for a normal crossing divisor $S_{Z'}$.
Then then the restriction of the sheaf $\sP$ to $S^{\mu(f)}(\Omega^1_Z(\log S_Z))$
lies in $S^{\mu(f)}(\pi^*(\Omega^1_{Z'}(\log S_{Z'})))$.
\end{enumerate}
\end{theorem}
\begin{proof} Parts i) and ii) have been shown in \cite{VZ3}, 1.4,
for canonically polarized manifolds with
$\mu(f)$ replaced by the fibre dimension $n$. We will just sketch the
changes which allow to extend the arguments used in \cite{VZ3}
to cover \ref{stthm2}, ii), iii) and iv), for canonically polarized manifolds.
Next we will try to convince the reader, that the same proof goes through
for minimal models of Kodaira dimension zero.

As in \cite{VZ3} we drop the assumption that $Y$ is projective. Leaving out a codimension
two subscheme, we may assume that $f$ is flat and that $\Delta$ is a relative normal
crossing divisor. Then we have the tautological exact sequence
$$
0 \>>> f^*\Omega^1_Y(\log S) \>>> \Omega^1_X(\log \Delta)
\>>> \Omega^1_{X/Y}(\log \Delta) \>>> 0
$$
and the wedge product sequences
\begin{multline}\label{exact}
0\>>> {f}^*\Omega^1_Y(\log S)\otimes
\Omega^{p-1}_{X/Y}(\log \Delta)
\>>> \\ {\mathfrak g \mathfrak r}(\Omega_X^p(\log
\Delta)) \>>> \Omega_{X/Y}^p(\log \Delta)\>>> 0,
\end{multline}
where
\begin{gather*}
{\mathfrak g \mathfrak r}(\Omega_X^p(\log \Delta))=
\Omega_X^p(\log \Delta)
/f^*\Omega^2_Y(\log S)\otimes \Omega^{p-2}_{X/Y}(\log
\Delta).
\end{gather*}
For the invertible sheaf
$\sL=\Omega_{X/Y}^n(\log \Delta)$ we consider the
sheaves
$$
F^{p,q}:=R^qf_*(\Omega^{p}_{X/Y}(\log
\Delta)\otimes\sL^{-1})
$$
together with the edge morphisms
$$
\tau_{p,q}:F^{p,q}\>>> F^{p-1,q+1}\otimes \Omega^1_{Y}(\log
S),
$$
induced by the exact sequence (\ref{exact}), tensored with
$\sL^{-1}$. As explained in \cite{VZ3}, Proof of 4.4 iii),
over $U$ the edge morphisms $\tau_{p,q}$ can also be obtained
in the following way. Consider the exact sequence
$$
0 \>>> T_{V/U} \>>> T_V \>>> f^*T_U \>>> 0,
$$
and the induced wedge product sequences
$$
0 \>>> T^{n-p+1}_{V/U} \>>> \tilde T^{n-p+1}_V \>>> T^{n-p}_{V/U}\otimes f^*T_U \>>> 0,
$$
where $\tilde{T}^{n-p+1}_{V}$ is a subsheaf of $T^{n-p+1}_V$. One finds edge morphisms
$$
\tau^\vee_{p,q}:(R^{q}f_*T^{n-p}_{V/U}) \otimes T_U \>>> R^{q+1}f_*T^{n-p+1}_{V/U}.
$$
Restricted to $\eta$ those are just the wedge product with the Kodaira-Spencer class.
Moreover, tensoring with $\Omega_U^1$ one gets back $\tau_{p,q}|_U$.
Hence $\mu(f)$ is the smallest number $m$ for which the composite
\begin{multline*}
\tau^m:F^{n,0}=\sO_Y \> \tau_{n,0} >> F^{n-1,1}\otimes \Omega_{U}^1
\> \tau_{n-1,1} >>
F^{n-2,2}\otimes S^2(\Omega^1_U) \>>> \cdots \\
\> \tau_{n-m+1,m-1}>> F^{n-m,m}\otimes S^m(\Omega^1_U)
\end{multline*}
is non-zero. Next we used that (replacing $Y$ by some covering)
there is an ample
invertible sheaf $\sA$ on $Y$ such that the kernel $\sK$ of
$$
{\rm id}_{\sA}\otimes\tau_{n-m,m}: \sA\otimes F^{n-m,m} \>>>
\sA\otimes F^{n-m-1,m+1} \otimes \Omega^1_{Y}(\log S)
$$
is negative, or to be more precise, that its dual is weakly positive.
This gives a non-trivial map
$$
\upsilon:\sA \otimes \sK^\vee \>>> S^m(\Omega^1_Y(\log S))
$$
and we take for $\sP$ its image.

The number $m$ was used in the proof of \cite{VZ3}, 1.4 ii),
hence there is no harm to replace the upper bound $n$, used there,
by the more precise number $\mu(f)$ in \ref{stthm2}, ii).

The sheaves $F^{p,q}$ are compatible with restriction to the subvariety
$Z$. The assumption $\mu(f)=\mu(h)$ implies that the restriction
$$
\upsilon|_Z:\sA|_{Z} \otimes \sK^\vee|_{Z} \>>> S^m(\Omega^1_{Z}(\log S_Z))
$$
is non-trivial. In fact, the kernel $\sK'$ of
$$
{\rm id}_{\sA_{Z}}\otimes\tau^{Z}_{n-m,m}: \sA|_{Z}\otimes F^{n-m,m}|_{Z}
\>>> \sA|_{Z} \otimes F^{n-m-1,m+1}|_{Z} \otimes \Omega^1_{Z}(\log S_Z)
$$
contains $\sK|_{Z}$, and the diagram
$$
\begin{CD}
\sA|_{Z} \otimes \sK^\vee|_{Z} \>>> S^m(\Omega^1_{Y}(\log S))|_{Z}\\
\A A A \V V V \\
\sA|_{Z} \otimes {\sK'}^\vee \>>> S^m(\Omega^1_{Z}(\log S_Z))
\end{CD}
$$
is commutative. One obtains iii).\\

Since the sheaves $F^{p,q}$ and the maps $\tau_{p,q}$ are
compatible with pullbacks, under the additional assumptions made
in iv), the image of
$$
\sA|_{Z} \otimes {\sK'}^\vee \>>>
S^m(\Omega^1_{Z}(\log S_Z))
$$
lies in $S^m(\pi^*(\Omega^1_{Z'}(\log S_{Z'})))$ and the same holds true
for the restriction of $\sP$.\\

If one considers the proof of \cite{VZ3}, 1.4, i) and ii), the assumption
that the fibres are canonically polarized is used twice. First of all, since we
apply in the proof of 4.4, iv), the Akizuki-Kodaira-Nakano vanishing theorem
to the restriction of $\omega_F$ to a smooth multicanonical divisor $B$.
If some power of $\omega_F$ is trivial the divisor $B$ is empty,
and there is nothing to show.

The second time is in the proof of \cite{VZ3}, 4.8. We use the diagram (2.8.1) and the fact
that the morphism $Z^{\#} \to Y^{\#}$ considered there is of maximal variation.
The construction of (2.8.1) just uses the existence of the moduli scheme
$M_h$, and it provides a morphism $Z^{\#} \to Y^{\#}$ induced by a
generically finite morphism $Y^{\#} \to M_h$.

This construction works in particular for the moduli scheme of polarized manifolds
$F$ with $\omega_F^\delta=\sO_F$, and \ref{stthm} implies that the variation of
the morphism $Z^{\#} \to Y^{\#}$ is again maximal.

The rest of the arguments, as given on page 311--313 of \cite{VZ3} remain unchanged,
and one obtains \ref{stthm2}, i), ii) and iii).
\end{proof}
For families $f:X\to Y$ with $\omega_{X_\eta}$ semiample,
and with $\mu(f)=n$ one can add to \cite{VZ3}, 1.4, a statement
similar to \ref{stthm2}, iii) and iv). Since the later will not be used, we omit it.
\begin{theorem}\label{stthm3} Assume
$\omega_{X_\eta}$ is semiample, and ${\rm Var}(f)=\dim(f)$.
\begin{enumerate}
\item[i.] There
exists a non-singular finite covering $\psi:Y'\to Y$ and
a big coherent subsheaf $\sP'$ of $\psi^*S^{m}(\Omega^1_Y(\log S))$,
for some $m \leq \mu(f)$.
\item[ii.] If $\mu(f)=n$, then one finds $m=n$ in i).
\item[iii.] Let $Z$ be a submanifold of $Y$ such that $S_Z=S\cap Z$
remains a normal crossing divisor, and such that
$W=X\times_YZ$ is non-singular. For the induced family
$h:W \to Z$ assume that $\mu(f)=\mu(h)=n$. Then
one can choose the covering $\psi$ such that $\psi^{-1}(Z)$
is non-singular, $\psi^{-1}(S_Z)$ a normal crossing divisor
and such that the image of the sheaf $\sP'$ from part i) in
$\psi^*S^n(\Omega^1_Z(\log S_Z))$ is non trivial.
\end{enumerate}
\end{theorem}
\begin{proof} We keep the notations from the sketch of the proof
of \ref{stthm2}.
For part i) we replaced in \cite{VZ3},
page 309 and 310, the sheaves $F^{p,q}$ (in fact a twist of those
by some invertible sheaf on $Y$) by some quotient sheaves. But then
$\mu(f)$ remains an upper bound for the number $m$,
used there, and one obtains \ref{stthm2}, i), as stated.

However, one has no control on the behavior of $m$ under restriction to
subvarieties. So for part ii) and iii) we have to recall the construction in more
detail. To get the weak positivity of the kernels $\sK$ one has to replace
(over some covering $Y'$ of $Y$ whose ramification divisor is in general
position) $\sA\otimes F^{n-m,m}$ by its image
$\sA\otimes \tilde F^{n-m,m}$  in some larger sheaf $E^{n-m,m}$.
Here
$$
\big(\bigoplus_{p+q=n}E^{p,q},\bigoplus_{p+q=n}\theta_{p,q}\big)
$$
is again a Higgs bundle, and $\theta_{p,q}$ is compatible with $\tau_{p,q}$.

As stated in the proof of \cite{VZ3}, 4.4, iv) the kernel and
cokernel of the map $$ \sA\otimes F^{n-m,m} \to E^{n-m,m} $$ are
direct images of the $n-m-1$-forms of a multicanonical divisor. If
$n=m$ there are no such forms, and $\sA\otimes F^{0,n}$ is a
subsheaf of $E^{0,n}$. Hence $\tau^n\neq 0$ implies that the
corresponding map for $\sA\otimes \tilde F^{n-m,m}$ is non-zero.
The compatibility with restrictions follows by the argument used
in the proof of \ref{stthm2}, iii).
\end{proof}
\section{Subvarieties of the moduli stack of
polarized manifolds of Kodaira dimension zero}\label{sectappl}

Theorem \ref{stthm2} has a number of geometric implication for manifolds
$U$ mapping to moduli stacks of polarized manifolds, i.e. for
morphisms $\varphi:U \to M_h$ induced by a family $f:V\to U$.
Those had been shown in \cite{VZ3} for the moduli stack
of canonically polarized manifolds. The proves are all based on
vanishing theorems for logarithmic differential forms, and they do not
refer to the type of fibres of $f$, once \ref{stthm2}, i) and ii),
is established.

Using \ref{stthm} we extended \ref{stthm2}, i) and ii), to a larger class
of families of polarized manifolds. Hence the geometric implications
carry over to this larger class, i.e. to the moduli stack of
polarized manifolds with $\omega_F^\delta=\sO_F$, for some $\delta>0$.
For the readers convenience we recall the statements below.

\begin{theorem}[see \cite{VZ3}, 5.2, 5.3, 7.2, 6.4, and 6.7]\label{cor1}
Let $M_h$ be the moduli scheme of canonically polarized manifolds,
or of polarized manifolds $F$ with $\omega_F^\delta$ trivial for some $\delta >0$.
\begin{enumerate}
\item[I.]
Assume that $U$ satisfies one of the following conditions
\begin{enumerate}
\item[a)] $U$ has a non-singular projective compactification $Y$
with $S=Y\setminus U$ a normal crossing divisor and with
boundary $T_Y(-\log S)$ weakly positive.
\item[b)] Let $H_1+ \cdots + H_\ell$ be a reduced normal crossing
divisor in $\bP^N$, and $\ell< \frac{N}{2}$. For $0\leq r\leq l$ define
\begin{gather*}
H = \bigcap_{j=r+1}^\ell H_j, \ \ \ S_i = H_i|_H, \ \ \
S = \sum_{i=1}^r S_i,
\end{gather*}
and assume $U= H \setminus S$.
\item[c)] $U=\bP^N\setminus S$ for a reduced normal crossing divisor
$$
S=S_1+ \cdots + S_\ell
$$
in $\bP^N$, with $\ell< N.$
\end{enumerate}
Then a morphism $U \to M_h$, induced by a family, must be
trivial.
\item[II.]
For $Y=\bP^{\nu_1}\times \cdots \times
\bP^{\nu_k}$ let
$$
D^{(\nu_i)}=D_0^{(\nu_i)}+\cdots +D_{\nu_i}^{(\nu_i)}
$$
be coordinate axes in $\bP^{\nu_i}$ and
$$
D=\sum_{i=1}^kD^{(\nu_i)}.
$$
Assume that $S=S_1+\cdots S_\ell$ is a divisor, such that $D+S$
is a reduced normal crossing divisor, and $\ell < \dim(Y)$.
Then there exists no morphism $\varphi:U=Y\setminus (D+S) \to M_h$
with
$$
\dim(\varphi(U)) > {\rm Max}\{\dim(Y)-\nu_i; \ i=1,\ldots ,k\}.
$$
\item[III.]
Let $U$ be a quasi-projective variety and let $\varphi:U \to
M_h$ be a quasi-finite morphism, induced by a family. Then
$U$ can not be isomorphic to the product of more than
$\mu(f)$ varieties of positive dimension.
\end{enumerate}
\end{theorem}

\section{Rigidity}\label{rigidity}
Again, $f:V\to U$ denotes a smooth family of manifolds with $\omega_{V/U}$
$f$-semiample and with ${\rm Var}(f)=\dim(U)>0$. We say that $f$ is rigid,
if there exists no non-trivial deformation over a non-singular
quasiprojective curve $T$.

Here a deformation of $f$ over
$T$, with $0\in T$ a base point, is a smooth projective morphism
$$g:\sV \to U\times T$$ for which there exists a commutative diagram
$$
\begin{CD}
V \> \simeq >> g^{-1}(U\times\{0\}) \>\subset >> \sV \\
\V f VV \V V V \V V g V \\
U \> \simeq >> U\times\{0\} \> \subset >> U\times T
\end{CD}.
$$
If the fibres $F$ of $f$ are
canonically polarized, or if some power of $\omega_F$ is trivial,
this says that morphisms from $U$ to the moduli stack do not deform.

\begin{proposition}\label{proprigid}
Assume either that $\omega_{X_\eta}$ is ample, or that
$\omega^\delta_{X_\eta}=\sO_{X_\eta}$, for some $\delta$.
Assume that ${\rm Var}(f)=\dim(U) >0$.
Let $T$ be a non-singular quasi-projective curve.
Let $g:\sV \to U\times T$ be a deformation of $f$.
If $\mu(f)=\mu(g)$, then ${\rm Var}(g)=\dim(U)$.
\end{proposition}
\begin{proof}
Suppose that ${\rm Var}(g)>\dim (U).$ Then
$$
\dim (U)+1 = \dim U\times T \geq  {\rm Var}(g)>\dim (U),
$$
hence, $\dim(U\times T)={\rm Var}(g).$ Let $\bar T$ be a non-singular
compactification of $T$, $S_{\bar T}=\bar T \setminus T$.
Correspondingly we write $S_{Y\times \bar T}$
for the complement of $U\times T$ in $Y\times \bar T$.

By Theorem \ref{stthm2}, ii, one finds a big coherent subsheaf
$\sP$ of
$$
S^{\mu(g)}(\Omega^1_{Y\times \bar T}(\log S_{Y\times \bar T})),
$$
and by \ref{stthm2}, iii) the image of $\sP$ in
$$
S^{\mu(g)}(\Omega^1_{Y\times \{0\}}(\log S_{Y\times \{0\}}))=
pr_1^*(S^{\mu(g)}(\Omega^1_{Y}(\log S)))|_{Y\times\{0\}}
$$
is non-zero. Then, the image of $\sP$ in
$$
pr_1^*(S^{\mu(g)}(\Omega^1_{Y}(\log S)))
$$
is non-zero, and for a point $y \in Y$ in general position, the image
of $\sP$ under
$$
S^{\mu(g)}(\Omega^1_{Y\times \bar T}
(\log S_{Y\times \bar T}))|_{\{y\} \times T}
\>>>
pr_1^*(S^{\mu(g)}(\Omega^1_Y(\log S)))|_{ \{y \}\times T}
$$
is not zero.
Note that any non-zero quotient of the coherent sheaf
$\sP|_{\{y\} \times T}$ for $y$ in general position must be big.
In fact, if $\sP$ is ample over some open dense subset $W_0$ of
$Y\times \bar T$, one just has to make sure that $\{y\}\times T$ meets $W_0$.
Since $pr_1^*(S^{\mu(g)}(\Omega^1_Y(\log S)))|_{\{y\}\times T}$ is a
direct sum of copies of $\sO_{\{y\} \times T}$ this is not possible.
\end{proof}
Using \ref{stthm3} instead of \ref{stthm2} one obtains a similar
result for families with $\omega_{V/U}$ $f$ semiample, whenever
$\mu(f)=n$.
\begin{proposition}\label{proprigidsa}
Assume that $\omega_{X_\eta}$ is semiample, that $${\rm
Var}(f)=\dim(U) >0$$
and that $$\mu(f)=\dim(X_\eta)=n.$$ Let $T$ be
a non-singular quasi-projective curve, and let $g:\sV \to U\times
T$ be a deformation of $f$. Then ${\rm Var}(g)=\dim(U)$.
\end{proposition}
\begin{proof}
If ${\rm Var}(g)>\dim (U)$, again one finds
$\dim(U\times T)={\rm Var}(g).$
Let us keep the notations from the proof of \ref{proprigid}.
By Theorem \ref{stthm3}, ii, one finds a finite covering
$\psi:Y'\to Y$ and a big coherent subsheaf
$\sP'$ of
$$
\psi^*(S^{\mu(g)}(\Omega^1_{Y\times \bar T}(\log S_{Y\times \bar T}))),
$$
and by \ref{stthm3}, iii) the image of $\sP'$ in
$$
\psi^*(S^{\mu(g)}(\Omega^1_{Y\times \{0\}}(\log S_{Y\times \{0\}})))=
\psi^*pr_1^*(S^{\mu(g)}(\Omega^1_{Y}(\log S))|_{Y\times\{0\}})
$$
is non-zero. Then, for a point $y \in Y$, in general position, the image
of $\sP'$ under
\begin{multline*}
\psi^*(S^{\mu(g)}(\Omega^1_{Y\times \bar T}
(\log S_{Y\times \bar T}))|_{\{y\} \times T})
\>>>\\
\psi^*(pr_1^*(S^{\mu(g)}(\Omega^1_{Y}(\log S)))|_{\{y\} \times T})
\end{multline*}
is not zero. Again, since the sheaf on the right hand side is trivial,
one obtains a contradiction.
\end{proof}

Let $h:\sX \to Z$ be a polarized family of manifolds $F$ with $\omega_F$
semiample and of maximal variation, over a non-singular quasi-projective
manifold $Z$. Assume that $\bar Z$ is a projective compactification of $Z$,
such that $\bar Z \setminus Z$ is a normal crossing divisor.
Assume in addition, that there is an open dense subscheme $Z_0$
such that for all subvarieties $\sU$ of $Z$ meeting $Z_0$
$$
{\rm Var}(\sX\times_{Z}\sU \to \sU)=\dim(\sU).
$$
Let $Y$ be a non-singular projective curve and let $U\subset Y$ be open
and dense. Let us write
$$
{\rm \bf H}={\rm\bf Hom}((Y,U),(\bar Z,Z ))
$$
for the scheme parameterizing non-trivial morphisms $\psi:Y \to \bar Z$ with
$\psi(U)\subset Z$ and
$$
{\rm \bf H}_{Z_0} ={\rm\bf Hom}((Y,U),(\bar Z,Z );Z_0) \subset {\rm\bf H}
$$
for those with $\psi(U)\cap Z_0 \neq \emptyset$.
Based on the bounds obtained in \cite{VZ1} we have shown in
\cite{VZ3} that ${\rm\bf H}_{Z_0}$ is of finite type.

\begin{corollary}\label{corrigid} \
\begin{enumerate}
\item[I.]
\begin{enumerate}
\item[a.] Let $\psi:U\to Z$ be a morphism and $f:V\to U$ the pull back
family. Assume that $\psi(U)\cap Z_0 \neq \emptyset$ and that
$$
\mu(f:V \to U)=\dim(F)=n
$$
Then the point $[\psi:Y\to \bar Z]$ is isolated in ${\rm\bf H}_{Z_0}$.
\item[b.] Assume for all fibres $F$ of $h^{-1}(Z_0) \to Z_0$ and for all
$\xi\in H^1(F,T_F)$
$$
0\neq \wedge^n \xi \in
H^n (F,\omega_F^{-1}).
$$
Then ${\rm\bf H}_{Z_0}$ is a finite set of points.
\item[II.] Assume that the fibres $F$ of $h:\sX\to Z$ are either
canonically polarized, or of Kodaira dimension zero.
\item[a.] Let $\psi:U\to Z$ be a morphism and $f:V\to U$ the pull back
family. Assume that $\psi(U)\cap Z_0 \neq \emptyset$ and that
$$
\mu(f:V \to U)=\mu(h:\sX \to Z).
$$
Then the point $[\psi:Y\to \bar Z]$ is isolated in ${\rm\bf H}_{Z_0}$.
\item[b.] Assume there exists a constant $\mu$ such
that for all fibres $F$ of $h^{-1}(Z_0) \to Z_0$ and for all
$\xi\in H^1(F,T_F)$
$$
0\neq \wedge^\mu \xi \in
H^\mu (F,T_F^{\mu})
$$
but
$$
0 = \wedge^{\mu+1} \xi \in
H^{\mu+1} (F,T_F^{\mu+1}).
$$
Then ${\rm\bf H}_{Z_0}$ is a finite set of points.
\end{enumerate}
\end{enumerate}
\end{corollary}
\begin{proof} In both cases
b) follows from a). For the latter assume that $[\psi]$
lies in a component of $\rm\bf H$ of dimension larger than zero.
Let $T$ be a curve in $\rm\bf H$, containing the point $[\psi]$. Then one has a
non-trivial deformation $\Psi:U\times T \to Z$ of $\psi$, hence a non-trivial
deformation $g:\sV \to U\times T$ of $f:V\to U$. By \ref{proprigid} in case II)
or by \ref{proprigidsa} in case I)
$$
{\rm Var}(g)={\rm Var}(f) < \dim(U\times T) = \dim(U)+1,
$$
contradicting the assumption made on $Z_0$ and $\sX \to Z$.
\end{proof}
Corollary \ref{corrigid}, II), should imply certain finiteness
results for curves in the moduli scheme $M_h$ of canonically
polarized manifolds, or the moduli scheme of minimal models of
Kodaira dimension zero meeting an open subscheme $W$ where the
assumption corresponding to the one in \ref{corrigid}, II), b),
holds true. However, one would have to show, that morphisms
$\varphi$ which factor through the moduli stack, are parameterized
by some coarse moduli scheme. Hopefully this can be done extending
the methods used in \cite{Cap} for moduli of curves to moduli of
higher dimensional manifolds.

Here we will show a slightly weaker statement, which coincides with
\ref{corrigid3} for $\mu=n$.

\begin{corollary}\label{corrigid2} \
Let $M_h$ be either the moduli scheme of canonically polarized
manifolds or the moduli scheme of polarized manifolds $F$ with
$\omega_F^\delta=\sO$ for some $\delta >0$. Let $0< \mu \leq
n=\dim(F)$ be a constant such that for all $(F,\sL)$, and for all
$\xi\in H^1(F,T_F)$ $$ 0 = \wedge^{\mu+1} \xi \in H^{\mu+1}
(F,T_F^{\mu+1}). $$ Then for a quasi-projective non-singular curve
$U$ there are only finitely many morphisms $\varphi:U\to M_h$
which are induced by a smooth family $f:V\to U$ with $\mu(f)=\mu$.
\end{corollary}
\begin{proof}
Let us choose any projective compactification $\bar M_h$ of $M_h$,
and an invertible sheaf $\sH$ on $\bar M_h$ which is ample with respect to $M_h$.
As usual, $Y$ will be a non-singular projective curve containing $U$.
We write $s$ for the number of points in $S=Y\setminus U$ and $g(Y)$ for the
genus of $Y$.

To show that there are only finitely many components of the scheme
$$ {\rm\bf Hom}((Y,U),(M_h,\bar M_h)) $$ which contain a morphism
$\varphi:U\to M_h$ factoring through the moduli stack, one has to
find an upper bound for $\varphi^*\sH$. To this aim one may assume
that $\bar M_h$ is reduced. The proof for the boundedness follows
the line of the proof of \cite{VZ3}, 6.2.

Kollar and Seshadri constructed (see \cite{Vie}, 9.25)
a finite covering of $M_h$ which factors through the moduli stack.

Consider any finite morphism $\pi:Z \to M_h$ with this property.
We choose a projective compactifications $\bar Z$ of $Z$ such that
$\pi$ extends to $\pi:\bar Z \to \bar{M}_h$. So $\pi^*\sH$ is
again ample with respect to $\pi^{-1}(M_h)$. Let $M_0$ be a
non-singular subvariety of $\pi(Z)\cap M_h$ with
$Z_0=\pi^{-1}(M_0)$ non singular.

Recall that for a family of projective varieties we constructed
in \cite{VZ3}, 2.7, a good open subset of the base space.
Applying this construction to the restriction of the universal family to $Z_0$,
we may assume furthermore, that $Z_0$ coincides with this subset.

By induction on the dimension of
$Z$, we may assume that we have found an upper bound for
$\varphi^*(\sH)$ whenever $\varphi(Y)\subset \pi(Z)\setminus M_0$.
Hence it is sufficient to find such a bound under the assumptions
that $\varphi(Y)\subset \pi(Z)$ and $\varphi(Y)\cap M_0 \neq
\emptyset$. There exists a finite covering $Y'$ of $Y$ of degree
$d \leq \deg(Z/\pi(Z))$, such that
$$
Y' \>\sigma >>  Y \> \varphi >> \pi(Z)
$$
factors through $\varphi':Y' \to Z$, and it is sufficient to bound
the degree of $\sigma^*\varphi^*\sH$. For simplicity, we assume that
$\varphi:Y \to \pi(Z)$ factors through $\varphi':Y\to Z$.

By \cite{VZ3}, 2.6 and 2.7, blowing up $Z$ with centers in
$Z\setminus Z_0$ we may assume that $Z$ is non singular, that
there exists a certain invertible sheaf $\lambda_\nu$ on $Z$, and a
constant $N_\nu >0$, such that
$$
\deg({\varphi'}^*\lambda_\nu) \leq
N_\nu\cdot \deg(\det(f_*\omega^\nu_{X/Y})),
$$
where, as usual, $f:X\to Y$ is
an extension of $V\to U$ to a projective manifold $X$. By the explicit description
of $\lambda_\nu$ in \cite{VZ3}, 2.6, d), and by \cite{VZ3}, 3.4, the sheaf
$\lambda_\nu$ is ample with respect to $Z_0$ for some $\nu >1$.
Hence it is sufficient to give
an upper bound for $\deg({\varphi'}^*\lambda_\nu)$, or for
$\deg(\det(f_*\omega^\nu_{X/Y}))$.

By \cite{VZ1} (see also \cite{BV} and \cite{Kov}) there exists a constant $e$,
depending only on the Hilbert polynomial $h$, with
\begin{equation*}
\deg(\det(f_*\omega_{X/Y}^\nu)) \leq
(n\cdot(2g(Y)-2+s)+s)\cdot\nu\cdot {\rm
rank}(f_*\omega_{X/Y}^\nu)\cdot e,
\end{equation*}
and we found the bound we were looking for.

It remains to show, that the points $$ [\varphi:Y\to \bar M_h] \in
{\rm\bf Hom}((Y,U),(M_h,\bar M_h)) $$ which are induced by a
family $f:X\to Y$ with $\mu(f)=\mu$, are discrete. If not, one
finds a positive dimensional manifold $T$ and a generically finite
morphism to ${\rm\bf Hom}((Y,U),(M_h,\bar M_h))$ whose image
contains a dense set of points where the corresponding morphism is
induced by a family. Let us choose a smooth projective
compactification $\bar T$ with $S_{\bar T} = \bar T \setminus T$ a
normal crossing divisor.

The induced morphism $Y\times T \to \bar M_h$ is not necessarily
factoring through the moduli stack, but using again the Koll\'ar
Seshadri construction again, we find a generically finite morphism
$\pi:Z \to Y\times \bar T$ which over $\pi^{-1}(U\times T)$ is
induced by a smooth family. Assume that $Z$ is non-singular and that
$S_Z=Z\setminus U\times T$ is a normal crossing divisor. Write
$p:Z \to \bar T$ for the induced morphism.

Applying \ref{stthm2}, ii), one obtains a big coherent subsheaf
$$\sP \subset S^\mu(\Omega^1_Z(\log S_Z)).$$ By part iii), is
image $\sP'$ in $S^\mu(\Omega_{Z/T}(\log S_Z))$ is non zero, and
iv) implies that for a dense set of points $t\in T$ the
restriction $\sP'|_{p^{-1}(t)}$ lies in $$
\pi^*S^\mu(\Omega^1_{Y\times \{t\}}(\log (S\times \{t\}))). $$

This is only possible, if $\sP'$ is a big subsheaf of $$
\pi^*S^\mu(\Omega^1_{Y\times \bar T}(\log (S\times \bar T))). $$
Restricting to $\pi^{-1}(\{y\}\times \bar T$, for general $y\in Y$
one obtains as in the proof of \ref{proprigidsa} a big subsheaf
of a trivial sheaf, a contradiction.
\end{proof}
Needless to say, Corollary \ref{corrigid2} is sort of empty, as long as
we do not know any answer to Problem \ref{problem2}.
\bibliographystyle{plain}

\end{document}